\documentclass{amsart}

\newcommand{\HW}{\mathrm{HW}}
\newcommand{\RND}{\mathrm{RND}}
\newcommand{\reg}{\text{reg}}

\newcommand{\End}{\operatorname{End}\nolimits}
\newcommand{\liea}[1]{\mathfrak{#1}}
\newcommand{\ad}{\operatorname{ad}\nolimits}
\newcommand{\rk}{\operatorname{rk}}
\newcommand{\lieg}[1]{\mathrm{#1}}
\newcommand{\im}{\operatorname{im}}
\newcommand{\CC}{{\mathbb C}}
\newcommand{\PP}{{\mathbb P}}
\newcommand{\Gr}{\operatorname{Gr}\nolimits}
\newcommand{\diag}{\operatorname{diag}\nolimits}
\newcommand{\cha}{\operatorname{char}}
\newcommand{\pafg}[2][]{\frac{\partial #1}{\partial #2}}
\newcommand{\ZZ}{{\mathbb Z}}

\theoremstyle{plain}
\newtheorem{thm}{Theorem}
\newtheorem{lm}[thm]{Lemma}

\newtheorem{prop}[thm]{Proposition}

\theoremstyle{definition}

\newtheorem{ex}[thm]{Example}

\begin{document}

\title{Small maximal spaces of non-invertible matrices}
\author{Jan Draisma}
\thanks{The author is supported by the Swiss National Science Foundation}
\address{Jan Draisma, Mathematisches Institut der Universit\"at Basel, 
	Switzerland}
\email{jan.draisma@unibas.ch}
\date{14 July 2004}

\maketitle

\section*{Abstract}

The {\em (generic) rank} of a vector space $\mathcal{A}$ of $n\times
n$-matrices is by definition the maximal rank of an element of
$\mathcal{A}$. The space $\mathcal{A}$ is called {\em rank-critical}
if any matrix space that properly contains $\mathcal{A}$ has a strictly
higher rank. I present a sufficient condition for rank-criticality,
and apply this condition to prove that the images of certain Lie algebra
representations are rank-critical.

A rather counter-intuitive consequence is that for infinitely many
$n$, there exists an $8$-dimensional rank-critical space of $n
\times n$-matrices having generic rank $n-1$, or, in other words: an
$8$-dimensional maximal space of non-invertible matrices. This settles
the question, posed by Fillmore, Laurie, and Radjavi in 1985, of whether
such a maximal space can have dimension smaller than $n$. As another
corollary, I prove that the image of the adjoint representation of a
semisimple Lie algebra is rank-critical.

\section{Results} \label{sec:Results}

This paper deals with linear subspaces of $\End(V)$, the space of
$K$-linear maps from an $n$-dimensional vector space $V$ over a field $K$
into itself. The {\em (generic) rank} of such a subspace $\mathcal{A}$,
denoted $\rk \mathcal{A}$, is by definition the highest rank of an
element of $\mathcal{A}$, and we call $\mathcal{A}$ {\em rank-critical}
if any linear subspace $\mathcal{B}$ of $\End(V)$ that properly contains
$\mathcal{A}$ has $\rk \mathcal{B} > \rk \mathcal{A}$. Note that a space
$\mathcal{A}$ is maximal among the {\em singular spaces}---that is, those
that only contain non-invertible matrices---if and only if $\mathcal{A}$
is rank-critical of rank $n-1$; in this case we call $\mathcal{A}$ a {\em
maximal singular space}. The main results of this paper are the following.

\begin{thm} \label{thm:Slm}
Let $K$ be a field of characteristic zero, let $m$ be an integer $\geq 3$,
and let $e$ be a positive integer. Then the image of the representation
of $\liea{sl}_m(K)$ on the space $V=K[x_1,\ldots,x_m]_{em}$ of homogeneous
polynomials of degree $em$ is a maximal singular subspace of $\End(V)$.
\end{thm}

In particular, taking $m=3$, we find that for every $n$ of the form
$\binom{3e+2}{2},\ e>1$ there exists an $8$-dimensional maximal singular
space of $n \times n$-matrices.

\begin{thm} \label{thm:Adjoint}
For any semisimple Lie algebra $\liea{g}$ over a field of characteristic
zero, $\ad \liea{g}$ is a rank-critical subspace of $\End(\liea{g})$.
\end{thm}

Theorems \ref{thm:Slm} and \ref{thm:Adjoint} are consequences of the
following proposition.

\begin{prop} \label{prop:SuffRND}
Let $\mathcal{A}$ be a subspace of $\End(V)$ and suppose that $|K| >
\rk\mathcal{A}=:r$. Set $\mathcal{A}_\reg:=\{A \in \mathcal{A} \mid \rk
A=r\}$ and define the space
\[ \RND(\mathcal{A}):=\{B \in \End(V) \mid B\ker A \subseteq \im A
	\text{ for all } A \in \mathcal{A}_\reg \}. \]
Then $\RND(\mathcal{A}) \supseteq \mathcal{A}$, and if equality holds,
then $\mathcal{A}$ is rank-critical. If, moreover, a group $G$ acts
linearly on $V$ and if $\mathcal{A}$ is stable under the conjugation
action of $G$ on $\End(V)$, then $\RND(\mathcal{A})$ is also a $G$-stable.
\end{prop}

The proof of Proposition \ref{prop:SuffRND} in Section
\ref{sec:Sufficient} bases on an elementary, but useful sufficient
condition for maximality of vector spaces in an {\em arbitrary} affine
variety embedded in a vector space. In Section \ref{sec:Construction},
we apply Proposition \ref{prop:SuffRND} to images of Lie algebra
representations; Theorems \ref{thm:Slm} and \ref{thm:Adjoint} are
proved there. Section \ref{sec:Random} lists some computer results
on rank-criticality of semisimple Lie algebra representations---in
particular, Theorem \ref{thm:Adjoint} came up as a conjecture using this
computer program.

\subsection*{Acknowledgments}
I thank Matthias B\"urgin, Arjeh Cohen, Hanspeter Kraft, Jochen Kuttler,
Martijn Stam, and Nolan Wallach for their help and for their interest
in the matter of this paper.

\section{Introduction and Motivation}

The direct motivation for this paper is the question, posed by Fillmore
{\em et al} in 1985 \cite{Fillmore85}, of whether a maximal singular
subspace of $\End(V)$ can have dimension smaller than $n$. I briefly
review three well-known constructions of maximal singular spaces that
led them to raise this question.

\begin{ex} \label{ex:UinW}
Fix subspaces $W,U$ of $V$ of dimensions $k-1$ and $k$,
respectively, and set $\mathcal{A}:=\{A \in \End(V) \mid AU
\subseteq W\}$. Then $\mathcal{A}$ is a singular space of dimension
$k(k-1)+(n-k)n=n^2-kn+k^2-k$. Moreover, it not hard to see that
$\mathcal{A}$ is maximal. We follow Eisenbud and Harris \cite{Eisenbud88}
in calling $\mathcal{A}$ and all its subspaces {\em compression spaces},
as they `compress' $U$ into $W$.
\end{ex}

\begin{ex} \label{ex:OddSkew}
Suppose that $n$ is odd, take $V=K^n$, and let $\mathcal{A}$ be the
space of all skew-symmetric matrices. As any skew-symmetric matrix has
even rank, the space $\mathcal{A}$ is singular, and it was observed in
\cite{Fillmore85} that $\mathcal{A}$ is maximal for all odd $n \geq
3$, under the assumption that $|K| \geq 3$. It is easy to see that
$\mathcal{A}$ is not a compression space.
\end{ex}

In both examples above, the dimension of $\mathcal{A}$ is quadratic in
$n$. An ingenious construction of smaller maximal singular spaces is
the following, attributed to Bob Par\'e in \cite{Fillmore85} and also
appearing in \cite{Lovasz89}.

\begin{ex} \label{ex:Skews}
Take $V=K^n$ and fix $n$ skew-symmetric $n \times n$-matrices
$A_1,\ldots,A_n$. Let $\phi$ be the linear map from $K^n$ into the space
$M_n(K)$ of $n \times n$-matrices over $K$ sending $x$ to the matrix
with columns $A_1x, A_2x, \ldots, A_nx$.  Then $\phi(K^n)$ is a singular
space in $M_n(K)$ because $x^t \phi(x)=0$ for all $x \in K^n$. Moreover,
in the particular case where $|K| \geq 3$, $A_i=E_{i,i+1}-E_{i+1,i}$
for $i<n$, and $A_n=E_{n,1}-E_{1,n}$, Fillmore {\em et al} showed that
$\mathcal{A}$ is maximal \cite{Fillmore85}.
\end{ex}

Many results in the literature exhibit sufficient conditions for a
singular space $\mathcal{A}$ to be a compression space: Dieudonn\'e
\cite{Dieudonne49} showed that every singular space of dimension $\geq
n^2-n$ either has a non-trivial common kernel or is dual to a space with
a common kernel. Under the assumption that $|K|$ is at least $2n-2$, this
result is sharpened as follows in \cite{Fillmore85}: if the dimension of
$\mathcal{A}$ is $>n^2-2n+2$ (which is the dimension of a compression
space with $k=2$), then $\mathcal{A}$ or its dual has a common kernel.
A condition of a different kind is that $\mathcal{A}$ be spanned by rank
one matrices; then a combinatorial argument shows that is a compression
space \cite{Lovasz89,Gelbord2002}. Analogues of these questions for
(skew-)symmetric matrices and for rank-critical spaces have also been
studied extensively in the literature \cite{Beasley87, Flanders62,
Gelbord2002, Ilic99, Loewy2001, Loewy94, Meshulam89}.

Yet another result of this type is part of the Kronecker-Weierstrass
theory of matrix pencils \cite{Gantmacher59} for $K=\CC$. I give a short
proof that is valid for other fields, as well.

\begin{prop}
Suppose that $|K| \geq n=\dim V$. Then any two-dimensional singular
subspace of $\End(V)$ is a compression space.
\end{prop}

\begin{proof}
If $A,B \in \End(V)$ are such that $KA+KB$ is singular, then the
homogeneous polynomial $p(s,t):=\det(sA+tB)$ vanishes on all $|K|+1$
points of $\PP^1(K)$. As $p(s,t)$ has degree $n<|K|+1$, it must be
identically zero, so that $sA+tB$ has a non-zero vector $u(s,t)$ in $K[s,t]
\otimes_K V$ in its kernel. But then any non-zero homogeneous component of
$u(s,t)$, say of degree $d$, is also annihilated by $sA+tB$; hence we
find $u_0,\ldots,u_d \in V$ such that $(sA+tB)(s^d u_0 + s^{d-1} t u_1 +
\ldots + t^d u_d)=0$, where we may assume that $u_0 \neq 0$. Taking the
of coefficients of $s^{d+1}, s^d t, \ldots, t^{d+1}$, we find
\[ Au_0=0,\ Au_1=-Bu_0,\ \ldots,\ Au_d=-Bu_{d-1},\text{ and } Bu_d=0.\] 
But then every element of $KA+KB$ maps the space $U:=\sum_i Ku_i$ into
the space $W:=\sum_i KAu_i$, which is strictly smaller because $A u_0=0$
while $u_0 \neq 0$.
\end{proof}

Example \ref{ex:2dim} below shows that the condition on $|K|$ in this
proposition is necessary.  This proposition, another proof of which is
given in \cite{Eisenbud88} for algebraically closed fields, shows that
if $|K| \geq n > 2$, then a maximal singular space cannot have dimension
$2$. On the other hand, Example \ref{ex:Skews} shows that there do exist
maximal singular spaces in $\End(V)$ of dimension $\dim(V)$, and this
led Fillmore {\em et al} to put forward their question above.

\section{Maximality of vector spaces in affine varieties} \label{sec:Sufficient}

Let $K$ be a field and let $M$ be a vector space over $K$. Denote by $L$
an algebraic closure of $K$ and set $M(L):=L \otimes_K M$. Let $Z$
be an affine algebraic variety in $M(L)$ defined over $K$, and let $N$
be a $K$-vector subspace of $M$ contained in $Z(K)$. We want a sufficient
condition for $N$ to be {\em maximal} among the subspaces of $M$ contained
in $Z(K)$. Therefore, set
\[ U:=\{m \in M(L) \mid N(L)+Lm \subseteq Z\}; \]
then $U$ is an affine variety defined over $K$. We make the following
assumption:
\begin{equation} \label{eq:Ass1}
\forall m \in M(L) : N + Km \subseteq Z(K) \Rightarrow m \in U(K); \tag{*}
\end{equation}
in particular, $0$ lies in $U(K)$ and so $N(L) \subseteq Z$.

\begin{lm} \label{lm:Deg}
Suppose that $Z$ is defined by polynomials $f_1,\ldots,f_k \in K[M]$.
Then $|K|>\max_i \deg f_i$ implies \eqref{eq:Ass1}. If, moreover, the
$f_i$ of highest degree are homogeneous, then $|K| \geq \max_i \deg f_i$
already implies \eqref{eq:Ass1}.
\end{lm}

\begin{proof}
Suppose that $m \in M$ satisfies $N+Km \subseteq Z(K)$, and let
$n_1,\ldots,n_d$ be a $K$-basis of $N$. Then we have for all $i$
\[ f_i(t_1n_1+\ldots+t_dn_d+sm) = 0 \text{ for all } t_1,\ldots,t_d,s
	\in K. \]
But a polynomial of degree $<|K|$ cannot vanish everywhere on an affine
space over $K$ unless it is the zero polynomial. Hence if $\deg f_i <|K|$,
then the left-hand side is the zero polynomial in $t_1,\ldots,t_d,s$
and vanishes on $L^{d+1}$, as well. We conclude that then $m \in U \cap
M=U(K)$, whence the first statement. For the second statement, note that
a {\em homogeneous} polynomial of degree $\leq K$ cannot vanish everywhere
on a projective space over $K$ unless it is the zero polynomial.
\end{proof}

By \eqref{eq:Ass1}, maximality of $N$ among the subspaces of $M$ contained
in $Z(K)$ is equivalent to $U(K)=N$, and a sufficient condition for this
is clearly $U=N(L)$. In principle, this condition can be verified using
Gr\"obner basis techniques, but $U$ may be hard to compute even for
moderately complicated $Z$. We therefore set out to find {\em linear}
sufficient conditions, as follows: Let $Z_\reg$ be the set of smooth
points of $Z$, and set $N_\reg := N \cap Z_\reg$, $N(L)_\reg := N(L)
\cap Z_\reg$. The second assumption we make is
\begin{equation} \label{eq:Ass2}
N_\reg \neq \emptyset. \tag{**}
\end{equation}
Now let 
\[ T_{N(L)_\reg} Z:=\bigcap_{n \in N(L)_\reg} T_n Z \]
be the intersection of all tangent spaces to $Z$ at points of $N(L)_\reg$,
where each $T_n Z$ is viewed as a vector subspace (through the origin)
of $M(L)$. Then we have the following lemma.

\begin{lm} \label{lm:SuffL}
The $L$-vector space $T_{N(L)_\reg} Z$ contains $U$. In particular, if
$T_{N(L)_\reg} Z=N(L)$, then $N(L)$ is maximal among the $L$-subspaces
of $M(L)$ contained in $Z$.
\end{lm}

\begin{proof}
For $m \in U$ and $n \in N(L)_\reg$ the line $\{n+tm \mid t \in L\}$
lies in $Z$, and {\em a fortiori} $m$ is tangent to $Z$ at $n$. We
conclude that $m \in T_{N(L)_\reg} Z$, as claimed. The second statement
is immediate.
\end{proof}

We want a $K$-rational version of this lemma. To this end, we define
\[ T_{N_\reg}(Z(K)):=\bigcap_{n \in N_\reg} 
	((T_n Z)(K))=(\bigcap_{n \in N_\reg} T_n Z)(K), \]
where the second equality follows from the fact that $Z$, and
therefore every $T_n Z$, is defined over $K$. Note that
$T_{N_\reg}(Z(K))$ contains $(T_{N(L)_\reg}Z)(K)$, as the intersection
is taken over a smaller set.

\begin{prop} \label{prop:SuffK}
The $K$-vector space $T_{N_\reg}(Z(K))$ contains $U(K)$. In particular, if
$T_{N_\reg}(Z(K))=N$, then $N$ is maximal among the $K$-vector subspaces
of $M$ contained in $Z(K)$.
\end{prop}

\begin{proof}
If $m \in U(K)$, then by the previous lemma $m \in (T_{N(L)_\reg}
Z)(K)$, which space is contained in $T_{N_\reg}(Z(K))$; this shows the
first statement. The second statement is now immediate from assumption
\eqref{eq:Ass1}.
\end{proof}

A randomised algorithm to compute the tangent space $T_{N_\reg} Z(K)$
bases on the following observation.

\begin{prop} \label{prop:Random}
For all non-negative integers $l$ and $e$, the set of
$(n_1,\ldots,n_l) \in N_\reg^l$ for which
\[ \dim_K \bigcap_{i=1}^l (T_{n_i} Z)(K) \geq e \]
is the set of $K$-rational points of a closed subvariety of $N(L)_\reg^l$.
\end{prop}

\begin{proof}
Let $d$ be the dimension of $Z$ and let $\gamma:Z_\reg \rightarrow
\Gr_d(M(L))$ be the Gauss map sending a point of $Z$ to its
tangent space. Now the set of all $l$-tuples $(T_1,\ldots,T_l)
\in \Gr_d(M(L))^l$ whose intersection has dimension at least $e$
is closed in $\Gr_d(M(L))^l$, hence so is its pre-image under
$(\gamma|_{N_\reg})^{\times l}$. The set of the proposition is the set
of $K$-rational points of this pre-image.
\end{proof}

For dimension reasons, the space $T_{N_\reg} Z(K)$ is the intersection
of {\em finitely many} tangent spaces $T_{n_i} Z(K),\ i=1,\ldots,m$
with $n_i \in N_\reg$. Now if $K$ is large (in particular, if $K$
is infinite), then the preceding proposition suggests the following
randomised algorithm to compute $T_{N_\reg} Z(K)$: First, find an upper
bound on $m$, and second, choose $m$ elements of $N$ at random. These
are probably smooth points by \eqref{eq:Ass2}, and by the preceding
proposition the intersection of their tangent spaces is probably equal to
$T_{N_\reg} Z(K)$. In particular, if this intersection is equal to $N$,
then you are {\em sure} that $N$ is a maximal vector space in $Z(K)$.
This algorithm was used to produce computational evidence for Theorem
\ref{thm:Adjoint}, see Section \ref{sec:Random}. For the second half of
Proposition \ref{prop:SuffRND} we need the following.

\begin{prop} \label{prop:Gmod}
Suppose that a group $G$ acts $K$-linearly on $M$, that $N$ is $G$-stable,
and that $Z$ is stable under the corresponding $L$-linear action on
$M(L)$. Then $T_{N_\reg} Z(K)$ is a $G$-submodule of $M$.
\end{prop}

\begin{proof}
Let $m \in T_{N_\reg} Z(K)$, $g \in G$, and $n \in N_\reg$. Then $g^{-1} n
\in N \cap Z_\reg=N_\reg$, and $g$ maps the $K$-rational tangent space
$(T_{g^{-1} n} Z)(K)$ isomorphically onto $(T_n Z)(K)$. As $m$ lies in
the former by assumption, $gm$ lies in the latter.  This shows that $gm
\in T_{N_\reg} Z(K)$.
\end{proof}

We will apply Proposition \ref{prop:SuffK} to the setting of Section
\ref{sec:Results}: $M=\End(V)$, where $V$ is a vector space over $K$,
$N=\mathcal{A}$ is a subspace of $\End(V)$ of generic rank $r$, and
$Z=R_r$ is the variety of $L$-linear maps $V(L) \rightarrow V(L)$ of rank
at most $r$. The following example shows that assumption \eqref{eq:Ass1}
is not automatic here.

\begin{ex} \label{ex:2dim}
Suppose that $K$ is a finite field with $q$ elements labelled
$c_1,\ldots,c_q$. Then the $2$-dimensional subspace $\mathcal{A}$ of
$M_{q+1}(K)$ spanned by the diagonal matrices $A=\diag(1,1,\ldots,1,0)$
and $B=\diag(c_1,c_2,\ldots,c_q,1)$ is singular---i.e., lies in
$R_q(K)$---but $\mathcal{A}(L)$ is not.
\end{ex}

We avoid this anomaly as follows: $R_r$ is defined by the $(r+1)
\times (r+1)$-minors of matrices, which are homogeneous polynomials of
degree $r+1$. Hence, if we assume that $|K| \geq r+1$, then assumption
\eqref{eq:Ass1} is automatically fulfilled by Lemma \ref{lm:Deg}; in
particular, the rank of $\mathcal{A}(L)$ is then also $r$.  Moreover,
the smooth points of $R_r$ are the linear maps of rank exactly $r$,
and by assumption $\mathcal{A}$ contains such maps, so that assumption
\eqref{eq:Ass2} is also satisfied. The space $T_{\mathcal{A}_\reg}
R_r(K)$ will be denoted $\RND(\mathcal{A})$, and its elements will be
called {\em rank-neutral directions} of $\mathcal{A}$.  We recall a
useful characterisation of the tangent spaces to $R_r$ at smooth points.

\begin{lm} \label{lm:Tangent}
For $A \in R_r(K)$ of rank $r$, the tangent space $(T_A R_r)(K)$ is
equal to $\{B \in \End(V) \mid B \ker A \subseteq \im A\}$.
\end{lm}

\begin{proof}
This is well known \cite[Example 14.16]{Harris92}; I only give an
intuitive argument: For $B$ to lie in $T_A R_r$, it is necessary and
sufficient that there be, for every $v \in \ker(A)$, a vector $w \in V$
for which $(A+\epsilon B)(v+\epsilon w)=0$ modulo $\epsilon^2$. The
coefficient of $\epsilon$ in this expression is $Aw+Bv$, so that the
existence of such a $w$ is equivalent to $Bv \in \im(A)$.
\end{proof}

Proposition \ref{prop:SuffRND} is now a direct consequence of
Proposition \ref{prop:SuffK}, Lemma \ref{lm:Tangent}, and Proposition
\ref{prop:Gmod}.

\begin{ex}
Proposition \ref{prop:SuffRND} provides another proof of rank-criticality
of compression spaces (under the condition $|K|\geq n$). Indeed, in
the notation of Example \ref{ex:UinW}, suppose that $B \in \End(V)$
does not map $U$ into $W$, and let $u \in U$ be such that $Bu \not \in
W$. It is not hard to construct an $A \in \mathcal{A}$ with $\ker A=Ku$
and $\im A \not \ni Bu$. But then $B \not \in (T_A R_{n-1})(K)$ by Lemma
\ref{lm:Tangent} and we conclude that $\RND(\mathcal{A})=\mathcal{A}$.
\end{ex}

\section{A construction of rank-critical spaces}
\label{sec:Construction}

The singular space of Example \ref{ex:OddSkew} is closed under
the commutator, and so are the compression spaces of Example
\ref{ex:UinW} if $W \subseteq U$.\footnote{Note that this property
is not preserved under the multiplication from the left and from
the right with arbitrary invertible matrices---however, deciding
whether a given subspace $\mathcal{A}$ of $\End(V)$ is $\lieg{GL}(V)
\times \lieg{GL}(V)$-conjugate to a Lie algebra is easily reduced to
a linear problem.} This suggests the study of the following situation:
Let $G$ be an affine algebraic group defined over $K$, and let $\rho:G
\rightarrow \lieg{GL}_K(V)$ be a finite-dimensional representation defined
over $K$. Let $\liea{g}$ be the Lie algebra of $G$, and denote the
corresponding representation $\liea{g} \rightarrow \End(V)$ also by
$\rho$. Set $r:=\rk \rho(\liea{g})$ and suppose that $|K| \geq r+1$.
Now $G$ acts on $\End(V)$ by $gA:=\rho(g)A\rho(g)^{-1}$, and both
$\rho(\liea{g})$ and the variety $R_r$ of linear maps of rank $\leq
r$ are $G$-stable. Proposition \ref{prop:SuffRND} implies: {\em
the rank-neutral directions of $\rho(\liea{g})$ form a $G$-module,
and if $\RND(\rho(\liea{g}))=\rho(\liea{g})$, then $\rho(\liea{g})$
is rank-critical.} In the rest of this section we assume that
$\cha K=0$, so that we can use the well-known representation theory of
semisimple Lie algebras to construct rank-critical spaces.

\begin{ex}
Recall Example \ref{ex:OddSkew}. Here $G$ is the group $\lieg{O}_n$
of orthogonal matrices, $\liea{g}=\liea{o}_n$, $\rho$ is the identity,
and $V=K^n$ is the standard $\liea{o}_n$-module. It is well known that
$\End(V)$ is the direct sum of three irreducible $\lieg{O}_n$-modules:
the space $\liea{o}_n$ of skew-symmetric matrices, the scalar multiples of
the identity $I$, and the space of symmetric matrices with trace $0$. We
now use Lemma \ref{lm:Tangent} to show that the last two modules are
not contained in $\RND(\liea{o}_n)$: Choose
\[ 
Y:=\diag(1,-1,0,\ldots,0) \text{ and }
X:=\left[ \begin{array}{c|c} 0 & 0 \\ \hline 0 & X' \end{array}
	\right],
\]
where $X' \in \liea{o}_{n-1}$ has full rank $n-1$. Then neither $I$
nor $Y$ maps $\ker X$ into $\im X$, hence the $\lieg{O}_n$-modules that
they represent are not contained in $\RND(\liea{o}_n)$. We conclude that
$\RND(\liea{o}_n)=\liea{o}_n$, so that $\liea{o}_n$ is maximal singular.
\end{ex}

Suppose now that $\liea{g}$ is semisimple and that it has a split Cartan
subalgebra $\liea{h}$. Then $V$ is the direct sum of its $\liea{h}$-weight
spaces $V_\lambda, \lambda \in \liea{h}^*$, and we have $r=\dim V-\dim
V_0$. Furthermore, we choose a Borel subalgebra $\liea{b}$ of $\liea{g}$
containing $\liea{h}$. We can then compute, for each $\liea{b}$-highest
weight $\lambda$ of $\End(V)$, the multiplicity of $\lambda$ among the
highest weights in $\RND(\rho(\liea{g}))$; see Section \ref{sec:Random}.

If, moreover, $K$ is algebraically closed, then the rank-neutral
directions of $\rho(\liea{g})$ can be characterised as follows:
\begin{equation} \label{eq:RNDLie}
\RND(\rho(\liea{g}))=\{Y \in \End(V) \mid (gY)V_0 \subseteq
	\bigoplus_{\lambda \neq 0} V_\lambda \text{ for all } g \in G\}. 
	\tag{***}
\end{equation}
Indeed, a generic element of $\rho(\liea{h})$ has kernel $V_0$ and image
$V_1:=\bigoplus_{\lambda \neq 0} V_\lambda$, so that the inclusion
$\subseteq$ follows from Lemma \ref{lm:Tangent}. On the other hand,
if any map in the $G$-orbit of $Y$ maps $V_0$ into $V_1$, then $Y$ is
tangent to $R_r$ at all points of $\{gX \mid g \in G, X \in \rho(\liea{h})
\text{ with } \ker X=V_0 \text{ and } \im X=V_1\}$.  This set is dense in
$\rho(\liea{g})_\reg$, and because the set $\{X \in \rho(\liea{g})_\reg
\mid Y \in T_X R_r\}$ is closed in $\rho(\liea{g})_\reg$, we conclude
that $Y \in \RND(\rho(\liea{g}))$.  We now proceed with the proof of
Theorem \ref{thm:Slm}.

Let $\liea{g}$ be the Lie algebra $\liea{sl}_m(K)$, let $k$ be a natural
number, and let $\rho:\liea{g} \rightarrow \End(V)$ be the representation
of $\liea{g}$ on the $k$-th symmetric power $V=S^k((K^m)^*)$ of
$(K^m)^*$, i.e., on the homogeneous polynomials on $K^m$ of degree
$k$. Let $\liea{h}\subseteq\liea{g}$ be the Cartan subalgebra of diagonal
matrices and let $\liea{b}\supseteq \liea{h}$ be the Borel subalgebra
of upper triangular matrices. The image of $\rho$ is spanned by the
(restrictions to $V$ of the) differential operators $x_i \pafg{x_j},
i \neq j$, and $x_i\pafg{x_i}-x_j \pafg{x_j}$, where the latter span
$\rho(\liea{h})$. The weight spaces in $V$ are one-dimensional and spanned
by the monomials $x_1^{a_1} \cdots x_m^{a_m}$ with $a_1+\ldots+a_m=k$.
The highest root vector in $\rho(\liea{g})$ is $x_m \pafg{x_1}$. To
apply Proposition \ref{prop:SuffRND}, we compute the highest weight
vectors in $\End(V)$.

\begin{lm} \label{lm:Symmetric}
The highest weight vectors in $\End(V)=V \otimes V^*$ are precisely
the powers $(x_m \pafg{x_1})^d$ for $d=0,\ldots,k$.
\end{lm}

\begin{proof}
It is clear that these are (non-zero) highest weight vectors; that there
are no others follows by a dimension computation. Alternatively, the lemma
is an easy application of the Littlewood-Richardson rule \cite{Fulton91}.
\end{proof}

The space $\rho(\liea{g})$ is singular if and only if $k$ is a multiple
of $m$, say $k=em$, and then $x_1^e x_2^e \ldots x_m^e$ spans the
zero weight space. The space $\rho(\liea{sl}_m)$ has no chance of being
maximal singular if $m=2$ (unless k=2), as then $\dim(V)=k+1=2e+1$ is odd
and $\lieg{SL}_2$ leaves invariant a non-degenerate symmetric bilinear
form on $V$; $\rho(\liea{g})$ is then contained in the larger singular
space of linear maps that are skew relative to this bilinear form.
This explains the condition $m\geq 3$ in Theorem \ref{thm:Slm}.

\begin{proof}[Proof of Theorem \ref{thm:Slm}]

First, if the image of $\liea{g}$ is maximal singular over a larger
field, then it is also maximal singular over $K$; so it suffices to
prove the theorem for $K$ algebraically closed. By Proposition
\ref{prop:SuffRND}, Lemma \ref{lm:Symmetric}, and the characterisation
\eqref{eq:RNDLie} of $\RND(\rho(\liea{g}))$ it suffices to prove that
if $d \neq 1$, then some element of the $\lieg{SL}_m$-orbit of $(x_m
\pafg{x_1})^d$ does {\em not} map $x_1^e \cdots x_m^e$ into the space
spanned by all {\em other} monomials.  This $\lieg{SL}_m$-orbit contains
the differential operators of the form $(x_1+x_2+x_3)^d(\alpha
\pafg{x_1}+\beta \pafg{x_2}+ \gamma \pafg{x_3})^d$ with
$\alpha+\beta+\gamma=0$, and the coefficient of $x_1^e \cdots x_m^e$ in
\[ (x_1+x_2+x_3)^d(\alpha \pafg{x_1}+\beta \pafg{x_2}+\gamma \pafg{x_3})^d 
	x_1^e \cdots x_m^e \]
is equal to
\[ 
\sum_{a,b,c:\ a+b+c=d} \left(\frac{d!}{a!\ b!\ c!}\right)^2 (e)_a (e)_b (e)_c\ 
	\alpha^a \beta^b \gamma^c,
\]
where $(e)_p:=e(e-1)\cdots(e-p+1)$ is the falling factorial. Lemma
\ref{lm:aplusbplusc} below shows that this polynomial in $\alpha,\beta$,
and $\gamma$ is not a multiple of $\alpha+\beta+\gamma$ if $d \neq
1$. In particular, for all $d \neq 1$ there exist $\alpha,\beta,\gamma
\in K$ with $\alpha+\beta+\gamma=0$ for which the coefficient above is
non-zero, that is, for which the term $x_1^e \cdots x_m^e$ does indeed occur in
$(x_1+x_2+x_3)^d(\alpha \pafg{x_1}+\beta \pafg{x_2}+\gamma
\pafg{x_3})^d x_1^e\cdots x_m^e$; this concludes the proof.
\end{proof}

\begin{lm} \label{lm:aplusbplusc}
For any integers $e>0$ and $d \in \{0\} \cup \{2,3,\ldots,3e\}$ the
polynomial
\[ P_{d,e}(\alpha,\beta,\gamma) := 
	\sum_{a,b,c:\ a+b+c=d} 
	\frac{d!}{a!\ b!\ c!} \binom{e}{a} \binom{e}{b} \binom{e}{c} 
	\alpha^a \beta^b \gamma^c \in \ZZ[\alpha,\beta,\gamma] 
\]
is {\em not} divisible by $\alpha+\beta+\gamma$.
\end{lm}

\begin{proof}
First, $P_{0,e}=1$ is not divisible by $\alpha+\beta+\gamma$; suppose next
that $d \geq 2e$. Then the highest monomial in $P_{d,e}$ with respect
to the lexicographic order with $\alpha>\beta>\gamma$ is $\alpha^e
\beta^e \gamma^{d-2e}$. Rewriting $P_{d,e}$ as a polynomial in the
elementary symmetric polynomials $\sigma_1:=\alpha+\beta+\gamma$,
$\sigma_2:=\alpha\beta+\beta\gamma+\gamma\alpha$, and
$\sigma_3:=\alpha\beta\gamma$ will therefore give the monomial
$\sigma_3^{d-2e}\sigma_2^{3e-d}$ a non-zero coefficient. Hence $P$
is not divisible by $\sigma_1$.

If $1< d < 2e$, then we claim that the coefficient of $\alpha^{d-1}
\beta$ in $P_{d,e}(\alpha,\beta,-\alpha-\beta)$ is non-zero. Indeed, this
coefficient is readily seen to equal $Q_{d-1}$, where 
\[ Q_d:=(e+1)(d+1)\sum_{a=0}^d \binom{d}{a} \binom{e}{a} \binom{e}{d-a}
	(-1)^a \frac{a}{a+1},\quad d \geq 1. \]
Automatic summation using Zeilberger-Wilf theory and the {\tt
Maple}-package {\tt EKHAD8} \cite{Brouwer2002, PWZ96} yields the following
expression for $Q_d$:
\[ Q_d=\begin{cases}
	(-1)^k \cfrac{(e+k)(e+k-1)\cdots(e-k+1)}{k\cdot((k-1)!)^2 }
		&\text{if }d=2k-1, \text{ and}\\
	(-1)^k (2e+1) \cfrac{(e+k)(e+k-1)\cdots(e-k+1)}{k\cdot((k-1)!)^2 }
		&\text{if }d=2k.
	\end{cases} \]
In particular, we find that $Q_d$ is non-zero for $1 \leq
k \leq e$, so that the coefficient of $\alpha^{d-1}\beta$ in
$P_{d,e}(\alpha,\beta,-\alpha-\beta)$ is in fact non-zero for $2\leq
d\leq 2e+1$. This concludes the proof of the lemma.
\end{proof}

We now prove rank-criticality of the images of adjoint
representations.

\begin{proof}[Proof of Theorem \ref{thm:Adjoint}.]
Again, it suffices to prove the theorem in the case where $K$ is
algebraically closed.  We will prove $\RND(\ad(\liea{g}))=\ad(\liea{g})$
and then apply Proposition \ref{prop:SuffRND}. Therefore, let $A$ be
a rank-neutral direction of $\ad(\liea{g})$ and let $x \in \liea{g}$
have centraliser $\liea{g}^x$ of minimal dimension; then $A \liea{g}^x
\subseteq [x,\liea{g}]$ by Lemma \ref{lm:Tangent}. In particular, we find
that if $[x,y]=0$, then the Killing form $\kappa$ of $\liea{g}$ vanishes
on $(x,Ay)$. As the commuting variety of $\liea{g}$ is irreducible
\cite{Richardson79}, this implies that $\kappa(x,Ay)=0$ for {\em all}
$x,y \in \liea{g}$ with $[x,y]=0$, independent of $\dim \liea{g}^x$. Now
consider the space
\[ M(\liea{g}):=\{ A \in \End(\liea{g}) \mid \kappa(x,Ay)=0
	\text{ for all } x,y \in \liea{g} \text{ with } [x,y]=0 \}; \]
we claim that it is equal to $\ad(\liea{g})$. First, assume that
this is true for {\em simple} $\liea{g}$, let $\liea{g}=\bigoplus_i
\liea{g}_i$ be a decomposition of $\liea{g}$ into simple ideals,
let $A \in M(\liea{g})$, and let $y \in \liea{g}_i$. Then $Ay$ is
$\kappa$-perpendicular to $\bigoplus_{j \neq i} \liea{g}_j$, so $Ay \in
\liea{g}_i$. In other words, every $\liea{g}_i$ is $A$-stable, and of
course $A|_{\liea{g}_i} \in M(\liea{g}_i)$. By assumption there exist
$z_i \in \liea{g}_i$ such that $A|_{\liea{g}_i}=\ad_{\liea{g}_i} z_i$,
and then $A=\sum_i \ad_\liea{g}(z_i)$.

It remains to prove $M(\liea{g})=\ad(\liea{g})$ for simple
$\liea{g}$. For $\liea{sl}_2$ this is easy, so we may suppose that
$\liea{g}$ has rank $\geq 2$. Setting $x=y$ in the condition on $A$, we
see that $M(\liea{g}) \subseteq \liea{o}(\kappa)$, the orthogonal Lie
algebra defined by $\kappa$. Moreover, $M(\liea{g})$ is stable under
conjugation with any automorphism of $\liea{g}$, and this implies two
things: first, that $M(\liea{g})$ is a $\liea{g}$-module and second,
using the Chevalley involution of $\liea{g}$, that it is self-dual as
such. Now the $\liea{g}$-module $\liea{o}(\kappa)/\ad(\liea{g})$ is
irreducible if $\liea{g}$ is not of type $A$---its highest weight can
be determined explicitly but is not of interest to us---while it is a
direct sum $W \oplus W^*$ for some non-self-dual module $W$ if $\liea{g}$
is of type $A$. In any case, $M(\liea{g})$ is either $\ad(\liea{g})$
or $\liea{o}(\kappa)$. But $M(\liea{g}) \neq \liea{o}(\kappa)$:
choose for instance $x,y$ in (a Cartan subalgebra of) $\liea{g}$
satisfying $[x,y]=0$ and $\kappa(x,x)=\kappa(y,y)=1-\kappa(x,y)=1$
and let $A$ be the map sending $x$ to $y$, $y$ to $-x$, and the
$\kappa$-orthogonal complement of $\langle x,y \rangle_K$ to $0$;
then $A \in \liea{o}(\kappa) \setminus M(\liea{g})$.  We conclude that
$M(\liea{g})=\RND(\ad(\liea{g}))=\ad(\liea{g})$, as claimed.
\end{proof}

\section{Some computer results} \label{sec:Random}

To compute the $G$-module $\RND(\rho(\liea{g}))$, where $\rho:\liea{g}
\rightarrow \End(V)$ is a finite-dimensional representation of a split
semisimple Lie algebra $\liea{g}$ over a field $K$ of characteristic $0$,
we use the following algorithm based on Proposition \ref{prop:Random}:

\begin{enumerate}

\item Compute the non-zero highest weight spaces $\HW_\mu$ of $\End(V)$
relative to a Borel subalgebra $\liea{b}$ and a split Cartan subalgebra
$\liea{h}$ contained in $\liea{b}$; this is elementary linear algebra.

\item For each of them, say $\HW_\mu$ of dimension $l$, choose $l$
random elements $X_1,\ldots,X_l \in \rho(\liea{g})$, verify that they have
maximal rank, and compute the subspace
\[ \{Y \in \HW_\mu \mid Y\ker X_i \subseteq \im X_i \text{ for all } 
	i=1,\ldots,l \}. \]
\item The dimension of this space is the multiplicity of $\mu$ among
the highest weights in $\RND(\rho(\liea{g}))$. Well, more precisely,
it is an upper bound to this multiplicity, which with high probability
is sharp.  
\end{enumerate}

If we find that the multiplicity of every highest weight
in $\RND(\rho(\liea{g}))$ is equal to its multiplicity in
$\rho(\liea{g})$, then $\rho(\liea{g})$ is rank-critical by Proposition
\ref{prop:SuffRND}. I list some examples, found by computer experiments
with this algorithm.

\begin{enumerate}
\item The images of the adjoint representations of split simple Lie
algebras of types $A_1,\ldots,A_4,B_2,\ldots,B_4,C_3,C_4$, and $G_2$
were proved rank-critical with this algorithm. This computational
evidence led to the formulation of Theorem \ref{thm:Adjoint}.
\label{it:a}

\item The image of the $26$-dimensional representation of $F_4$ is
rank-critical (of rank $24$).

\item Let $\liea{g}$ be split simple of type $G_2$ and let $\rho$
be the $7$-dimensional representation of highest weight $[1,0]$, with
zero weight multiplicity $1$. Then $\RND(\rho(\liea{g}))$ is equal to
$\liea{o}_7$, which of course is still singular, so $\rho(\liea{g})$
is not a maximal singular space.

\item Similarly, if $\liea{g}$ is of type $G_2$ and $\rho$ is the
$27$-dimensional representation of highest weight $[2,0]$ with zero
weight multiplicity $3$, then $\RND(\rho(\liea{g}))$ is equal to
$\phi(\liea{o}_7)$, where $\phi$ is the representation of $\liea{o}_7$
of highest weight $[2,0,0]$, which restricts to $\rho$ on $\liea{g}$.
As both $\rho(\liea{g})$ and $\phi(\liea{o}_7)$ have generic rank $24$,
the former is not rank-critical but the latter is.

\item Let $\rho$ be the 35-dimensional irreducible representation of
$\liea{g}=\liea{sl}_3$ of highest weight $[4,1]$. Then $\RND(\rho(\liea{g}))$
is a sum of three irreducible modules of highest weights $[1,4],[1,1],$
and $[4,1]$. Hence, Proposition \ref{prop:SuffRND} cannot be applied
to conclude rank-criticality. Note that by the results of Dynkin
\cite{Dynkin57a} the image of $\rho$ is a maximal subalgebra of
$\liea{sl}_{35}$, so that there is no easy argument as in the previous
two examples, that shows that $\rho(\liea{g})$ is {\em not} rank-critical.
\end{enumerate}

Theorem \ref{thm:Adjoint} was first conjectured on the basis of
\eqref{it:a}.

\section{Conclusion and further questions} \label{sec:Conclusion}

Representations of semisimple Lie algebras yield an abundance
of rank-critical matrix spaces, which suggests that the (commonly
believed to be intractable) classification of such spaces may in some
sense include the classification of Lie algebra representations. In
particular, among the representations of $\liea{sl}_3$ we found infinitely
many where the image of $\liea{sl}_3$ is a maximal singular space.
The matrix spaces $\mathcal{A}$ constructed this way actually satisfy
an {\em a priori} stronger condition than rank-criticality, namely:
$\mathcal{A}=\RND(\mathcal{A})$. These results pose many questions for
further research, of which the following seem most interesting.

\begin{enumerate}

\item Directly describe, given the highest weights in a representation
$\rho$ of a split semisimple Lie algebra $\liea{g}$, the highest weights
in the $\liea{g}$-module $\RND(\rho(\liea{g}))$. The present proofs of
Theorems \ref{thm:Slm} and \ref{thm:Adjoint} are somewhat {\em ad hoc},
and the algorithm of Section \ref{sec:Random} is computationally rather
intensive, so that it only works for representations of dimensions at most
$50$ or so.

\item Investigate the discrepancy between rank-criticality and the
condition $\mathcal{A}=\RND(\mathcal{A})$.

\item For $m_n$ the minimal dimension of a maximal singular space of
$n \times n$-matrices, determine $\liminf_{n \rightarrow \infty} m_n$
(which is larger than $2$ and at most $8$, as we have seen, while---to
the best of my knowledge---it was previously believed to be infinite)
and $\limsup_{n \rightarrow \infty} m_n$.

\item Investigate whether the maximal singular spaces of Theorem
\ref{thm:Slm} remain maximal modulo primes, and, more generally, whether
rank-critical spaces constructed in Section \ref{sec:Construction}
remain rank-critical modulo primes.

\end{enumerate}


\begin{thebibliography}{10}

\bibitem{Beasley87}
LeRoy~B. Beasley.
\newblock {Null spaces of spaces of matrices of bounded rank.}
\newblock In {\em {Current trends in matrix theory, Proc. 3rd Conf.,
  Auburn/Ala. 1986}}, pages 45--50. North-Holland, New York, 1987.

\bibitem{Brouwer2002}
Andries Brouwer.
\newblock Automatic summation using {Z}eilberger-{W}ilf theory.
\newblock {\em Nieuw Arch. Wiskd. (5)}, 3(4):308--312, 2002.

\bibitem{Dieudonne49}
Jean Dieudonn\'e.
\newblock {Sur une g\'en\'eralisation du groupe orthogonal \`a quatre
  variables.}
\newblock {\em Arch. Math., Oberwolfach}, 1:282--287, 1949.

\bibitem{Dynkin57a}
Eugene~B. Dynkin.
\newblock Maximal subgroups of the classical groups.
\newblock {\em Am. Math. Soc. Transl.}, II. Ser. 6:245--378, 1957.

\bibitem{Eisenbud88}
David Eisenbud and Joe Harris.
\newblock {Vector spaces of matrices of low rank.}
\newblock {\em Adv. Math.}, 70(2):135--155, 1988.

\bibitem{Fillmore85}
P.~Fillmore, C.~Laurie, and H.~Radjavi.
\newblock On matrix spaces with zero determinant.
\newblock {\em Linear Multilinear Algebra}, 18:255--266, 1985.

\bibitem{Flanders62}
H.~Flanders.
\newblock {On spaces of linear transformations with bounded rank}.
\newblock {\em J. Lond. Math. Soc.}, 37:10--16, 1962.

\bibitem{Fulton91}
William Fulton and Joe Harris.
\newblock {\em Representation theory. {A} first course}.
\newblock Number 129 in Graduate Texts in Mathematics. Springer-Verlag, New
  York, 1991.

\bibitem{Gantmacher59}
F.R. Gantmacher.
\newblock {\em {The theory of matrices. Vol. 2.}}
\newblock {AMS Chelsea Publishing}, {New York}, 1959.

\bibitem{Gelbord2002}
Boaz Gelbord and Roy Meshulam.
\newblock {Spaces of singular matrices and matroid parity.}
\newblock {\em Eur. J. Comb.}, 23(4):389--397, 2002.

\bibitem{Harris92}
Joe Harris.
\newblock {\em {Algebraic geometry. A first course.}}
\newblock Number {133} in {Graduate Texts in Mathematics}. {Springer-Verlag},
  {Berlin etc.}, 1992.

\bibitem{Ilic99}
Bo~Ilic and J.M. Landsberg.
\newblock {On symmetric degeneracy loci, spaces of symmetric matrices of
  constant rank and dual varieties.}
\newblock {\em Math. Ann.}, 314(1):159--174, 1999.

\bibitem{Loewy2001}
Raphael Loewy.
\newblock {Large spaces of symmetric matrices of bounded rank are
  decomposable.}
\newblock {\em Linear Multilinear Algebra}, 48(4):355--382, 2001.

\bibitem{Loewy94}
Raphael Loewy and Nizar Radwan.
\newblock {Spaces of symmetric matrices of bounded rank.}
\newblock {\em {Linear Algebra Appl.}}, 197-198:189--215, 1994.

\bibitem{Lovasz89}
L\'aszl\'o Lov\'asz.
\newblock {Singular spaces of matrices and their application in combinatorics.}
\newblock {\em {Bol. Soc. Bras. Mat., Nova S\'er.}}, 20(1):87--99, 1989.

\bibitem{Meshulam89}
Roy Meshulam.
\newblock {On two extremal matrix problems.}
\newblock {\em {Linear Algebra Appl.}}, 114-115:261--271, 1989.

\bibitem{PWZ96}
Marko Petkov{\v s}ek, Herbert~S. Wilf, and Doron Zeilberger.
\newblock {\em $A=B$}.
\newblock A. K. Peters, Wellesley, Massachusetts, 1996.

\bibitem{Richardson79}
R.~W. Richardson.
\newblock Commuting varieties of semisimple lie algebras and algebraic groups.
\newblock {\em Compositio Math.}, 38(3):311--327, 1979.

\end{thebibliography}

\end{document}